\documentclass[leqno]{article}

\usepackage{graphicx}

\usepackage{latexsym}
\usepackage{amsfonts}
\usepackage[psamsfonts]{amssymb}
\usepackage{enumerate}

\RequirePackage[OT1]{fontenc}
\RequirePackage[amsthm,amsmath]{imsart}

\theoremstyle{plain} 
\newtheorem{theorem}{Theorem}

\newtheorem{lemma}{Lemma}
\newtheorem{proposition}
{Proposition}

\theoremstyle{definition} 

\theoremstyle{definition} 

\theoremstyle{remark} 

\theoremstyle{remark} 
\newtheorem{remark}
{Remark}
\newtheorem*{remark*}{Remark}

\newcommand{\esssup}{\operatorname{ess\,sup}}

\newcommand{\xx}{\mathbf{x}}

\renewcommand{\u}{\nearrow}
\renewcommand{\d}{\searrow}

\newcommand{\lp}{\left(}
\newcommand{\rp}{\right)}

\renewcommand{\le}{\leqslant}
\renewcommand{\ge}{\geqslant}

\newcommand{\la}{\lambda}

\newcommand{\vpi}{\varphi}
\renewcommand{\Psi}{\overline{\Phi}}

\renewcommand{\P}{\operatorname{\mathsf{P}}} 

\newcommand{\E}{\operatorname{\mathsf{E}}}

\newcommand{\R}{\mathbb{R}}

\newcommand{\vp}{\varepsilon}

\begin{document}


\begin{frontmatter}

\title{On inequalities for sums of bounded random variables}
\runtitle{Inequalities for bounded random variables}

\begin{aug}
\author{\fnms{Iosif} \snm{Pinelis}\ead[label=e1]{ipinelis@math.mtu.edu}}
\runauthor{Iosif Pinelis}


\address{Department of Mathematical Sciences\\
Michigan Technological University\\
Houghton, Michigan 49931, USA\\
E-mail: \printead[ipinelis@mtu.edu]{e1}}
\end{aug}






\begin{abstract}
Let $\eta_{1},\eta_2,\ldots$ be independent (not necessarily identically distributed) zero-mean random variables (r.v.'s) such
that $|\eta_i|\le1$ almost surely for all $i$, and let
$Z$ stand for a standard normal r.v.
Let $a_1,a_2,\dots$ be any real numbers such that
$a_1^2+a_2^2+\dots=1.$
It is shown that then 
$$
\P(a_1\eta_1+a_2\eta_2+\dots\ge x)  
\le \P(Z\ge x-\la/x) \quad\forall x>0,
$$
where $\la := \ln\frac{2e^3}9=1.495\dots$.
The proof relies on (i) another probability inequality and (ii) a l'Hospital-type rule for monotonicity, both developed elsewhere.
A multidimensional analogue of this result is given, based on a dimensionality reduction device, also developed elsewhere. 
In addition, extensions to (super)martin\-gales are indicated.
\end{abstract}







\end{frontmatter}



To begin with, let $\vp_1,\vp_2,\dots$ be independent Rademacher random variables (r.v.'s), so that $\P(\vp_i=1)=\P(\vp_i=-1)=\frac12$ for all $i$. Let $a_1,a_2,\dots$ be any real numbers such that
\begin{equation}
	a_1^2+a_2^2+\dots=1.\label{eq:a}
\end{equation}

Using a result due to Eaton~\cite{eaton}, Edelman \cite{edel} proposed an interesting inequality for normalized Rademacher sums:
\begin{equation}\label{eq:edel}
\P(a_1\vp_1+a_2\vp_2+\dots\ge x) \le \P\lp Z\ge x-1.5/x\rp\quad\text{for all $x>0$},
\end{equation}
where $Z$ is a standard normal r.v.
Employing certain conditioning, Edelman \cite{edel} also gave applications of inequality \eqref{eq:edel} to statistical inference based on Student's $t$ statistic. Before that, the same conditioning idea (in relation with an inequality due to Hoeffding~\cite{hoeff} in place of \eqref{eq:edel}) was given by Efron~\cite{efron} and then by Eaton and Efron~\cite{eaton-efron}, in more general settings. 

The sketch of proof offered in \cite{edel} for inequality \eqref{eq:edel} required an apparently nontrivial iterative computation procedure, which I have not been able to 
reproduce within a reasonable amount of computer time, because of rapid deterioration of precision at every step of the iterative procedure. 

In this note, inequality \eqref{eq:edel} is proved with a slightly better constant
$$\la := \ln\frac{2e^3}9=1.495\dots$$
in place of $1.5$. 
In fact, somewhat more general results will be proved here. 

Let $\eta_1,\eta_2,\ldots$ be independent (not necessarily identically distributed) zero-mean r.v.'s such
that $|\eta_i|\le1$ almost surely (a.s.) for all $i$.
It is assumed throughout that normalization condition \eqref{eq:a} holds.
Let 
$$S:=a_1\eta_1+a_2\eta_2+\dots.$$

\begin{theorem}\label{th:pin-edel} 
For all $x>0$,
\begin{align}
\P(S\ge x)  
&\le W(x):=\min\big(e^{-x^2/2},\P(Z\ge x-\la/x) \big);
\label{eq:pin-edel}\\
\P(|S|\ge x)  
&\le\widetilde W(x):=
\min\big(\tfrac1{x^2},\P(|Z|\ge x-\la/x) \big).
\label{eq:pin-edel tilde}
\end{align}
\end{theorem}

\begin{remark}
The upper bounds in \eqref{eq:pin-edel} and \eqref{eq:pin-edel tilde} hold e.g. for $\P(S_n\ge x)$ and $\P(|S_n|\ge x)$, respectively, $\forall n$, where $(S_i)$ is a martingale with $S_0=0$ a.s.\ and differences $X_i:=S_i-S_{i-1}$ ($i\ge1$) such that $\sum_i\esssup|X_i|^2\le1$. Other extensions hold as well; look in \cite{asymm} for appearances of the constant $c_{3,0}=2e^3/9=e^\la$ together with $\P(Z\ge\ldots)$ or $\P(|Z|\ge\ldots)$.
Using the dimensionality reduction device given in \cite{dim-reduct}, one can also obtain a multi-dimensional generalization of \eqref{eq:pin-edel tilde}:
\begin{equation*}
\P\lp  \|\eta_1\xx_1+\eta_2\xx_2+\dots\| \ge x \rp
\le\min\big(\tfrac1{x^2},\P(|Z|\ge x-\la/x) \big)\quad
\forall x>0, 
\end{equation*}
where 
$\xx_1,\xx_2,\dots$ are any non-random vectors in a Hilbert space $(H,\|\cdot\|)$ such that $\|\xx_1\|^2+\|\xx_2\|^2+\dots=1$. 
Cf. \cite[Remark~1.4]{asymm}. 
\end{remark}

Two-tail inequality \eqref{eq:pin-edel tilde} can be easily deduced from \eqref{eq:pin-edel}. Indeed, inequality $\P(|S|\ge x)\le\tfrac1{x^2}$ for $x>0$ follows from Markov's inequality, since $\E S ^2\le1$. 
As for inequality $\P(|S|\ge x)\le\P(|Z|\ge x-\la/x)$, it is trivial for $x\in(0,\sqrt\la]$, while for $x>\sqrt\la$ it obviously follows from \eqref{eq:pin-edel}. 

Therefore, to complete the proof of Theorem~\ref{th:pin-edel} it remains to prove inequality
\eqref{eq:pin-edel}, which is an immediate corollary of the following two propositions.

\begin{proposition}\label{prop:pin-eaton} 
For all $x>0$,
\begin{align*}
\P(S\ge x)  
&\le V(x):=\min\big(e^{-x^2/2},e^\la\P(Z\ge x) \big).
\end{align*}
\end{proposition}

\begin{proposition}\label{prop:less} 
For all $x>0$,
\begin{align*}
V(x) &\le W(x). 
\end{align*}
\end{proposition}

Proposition~\ref{prop:pin-eaton} is well known. Inequality $\P(S\ge x)\le e^{-x^2/2}$ for $x>0$ follows from a result due to Hoeffding~\cite{hoeff} and later improved in \cite{normal}. 
As for inequality 
\begin{align}
\P(S\ge x)&\le e^\la\P(Z\ge x)\quad\forall x\in\R,
\label{eq:pin}	
\end{align} 
its two-tail version, $\P(|S|\ge x)\le e^\la\P(|Z|\ge x)$ $\forall x\in\R$,
was given first in \cite{pin91,pin94}. The right-tail inequality \eqref{eq:pin} can be proved quite similarly; alternatively, it
follows from general results of \cite{pin98}. Recently, different generalizations of \eqref{eq:pin} were given e.g.\ in \cite[(1.16)]{bent-ap} and \cite[Corollary~3.4]{asymm}.

As for Proposition~\ref{prop:less}, it 
is an immediate corollary of the following three lemmas.

\begin{lemma}\label{lem:V} 
One has
	\begin{align}
V(x)&=\begin{cases}
e^{-x^2/2} &\text{if}\quad  0< x\le z_V,\\
e^\la\P(Z\ge x) &\text{if}\quad  z_V\le x<\infty,
\end{cases}
\label{eq:V}
\end{align}
where 
$$z_V=1.312\ldots$$
is the unique root of the equation $e^{-z^2/2}=e^\la\P(Z\ge z)$ for $z>0$.
\end{lemma}

\begin{lemma}\label{lem:W} 
One has
	\begin{equation}\label{eq:W}
W(x)=\begin{cases}
e^{-x^2/2} &\text{if}\quad  0\le x\le z_W,\\
\P(Z\ge x-\la/x) &\text{if}\quad  z_W\le x<\infty,
\end{cases}
\end{equation}
where 
$$z_W=1.365\dots$$
is the unique root of the equation $e^{-z^2/2}=\P(Z\ge z-\la/z)$ for $z>0$, so that $z_W>z_V$.
\end{lemma}

\begin{lemma}\label{lem:less} 
One has
\begin{align*}
	e^\la\P(Z\ge x)&\le\P(Z\ge x-\la/x)\quad\forall x\ge z_V.
\end{align*}
\end{lemma}

Concerning the upper bound $\widetilde W(x)$ in \eqref{eq:pin-edel tilde}, the following analogue of Lemma~\ref{lem:W} may be of interest.

\begin{proposition}\label{prop:W tilde}
One has 
\begin{equation*}
\widetilde W(x)=\begin{cases}
1 &\text{if}\quad 0<x\le1,\\
\tfrac1{x^2} &\text{if}\quad  1\le x\le z_{\widetilde W},\\
\P(|Z|\ge x-\la/x) &\text{if}\quad  z_{\widetilde W}\le x<\infty,
\end{cases}
\end{equation*}
where 
$$z_{\widetilde W}=1.865\dots$$
is the unique root of the equation $\tfrac1{z^2}=\P(|Z|\ge z-\la/z)$ for $z\in(\sqrt\la,\infty)$.
\end{proposition}

In turn, the proofs of Lemmas~\ref{lem:V}, \ref{lem:W}, \ref{lem:less} and Proposition~\ref{prop:W tilde}
rely on the following particular cases of the l'Hospital-type rules for monotonicity given in \cite[Propositions~4.1 and 4.3; see also Remark~5.5]{borwein}. 

\begin{proposition}\label{prop:special-ud} 
\emph{\cite{borwein}}\quad
Let $-\infty\le a<b\le\infty$. Let $f$ and $g$ be real-valued 
differentiable functions,
defined on the interval $(a,b)$,
such that $g$ and $g'$ do not take on the zero value on $(a,b)$.
Let 
$$r:=f/g\quad\text{and}\quad\rho:=f'/g'.$$   
Suppose that $f(b-)=g(b-)=0$. Then the following two statements are true.
\begin{description}
	\item[(i)] Suppose that $\rho\d$ (that is, $\rho$ is decreasing) on $(a,b)$. Then $r\d$ on $(a,b)$. 
		\item[(ii)] 
	Suppose that $\rho\u\d$ on $(a,b)$ -- that is, for some $c\in(a,b)$, $\rho\u$ ($\rho$ is increasing) on $(a,c)$ and $\rho\d$ on $(c,b)$. Then $r\d$ or $\u\d$ on $(a,b)$. 
\end{description}
\end{proposition}

Alternatively, one can use \cite[Theorem 1.16]{waves} instead of Proposition~\ref{prop:special-ud}.

\begin{proof}[Proof of Lemma~\ref{lem:V}]
Consider the ratios $r:=f/g$ and $\rho:=f'/g'$, where $f(x):=e^\la\P(Z\ge x)$ and $g(x):=e^{-x^2/2}$. 
It is apparently a well known fact that $r$ is decreasing on $(0,\infty)$. However, this follows immediately from part (i) of Proposition~\ref{prop:special-ud}, because $\rho(x)=\frac{e^\la}{x\sqrt{2\pi}}$ is decreasing on $(0,\infty)$. 
Besides, $r(0)=e^\la/2>1$ and, by l'Hospital's rule for limits, $r(\infty-)=\rho(\infty-)=0$. Now Lemma~\ref{lem:V} follows. 
\end{proof}

\begin{proof}[Proof of Lemma~\ref{lem:W}]
Consider the ratios $r:=f/g$ and $\rho:=f'/g'$, where $f(x):=\P(Z\ge x-\la/x)$ and $g(x):=e^{-x^2/2}$. 
Then
$$\rho(x)=\frac{\la+x^2}
{\sqrt{2\pi }\, x^3\, e^{\frac{\la ^2}{2 x^2}-\la }}
\quad\text{and}\quad
\rho '(x) 
   =\frac{\la ^3-(3-\la ) \la x^2-x^4}
   {\sqrt{2\pi }\, x^6\, e^{\frac{\la ^2}{2 x^2}-\la }},$$ 
so that $\rho'$ changes sign from $+$ to $-$ on $(0,\infty)$ and hence $\rho\u\d$ on $(0,\infty)$. 
Now, by part (ii) of Proposition~\ref{prop:special-ud}, $r\d$ or $\u\d$ on $(0,\infty)$. But $r(0+)=1$, $r(1)=1.13\ldots>1$, and, by l'Hospital's rule for limits, $r(\infty-)=\rho(\infty-)=0$. 
Therefore, $r\u\d$ on $(0,\infty)$. 
Now Lemma~\ref{lem:W} follows.     
\end{proof}

\begin{proof}[Proof of Lemma~\ref{lem:less}]
Consider the ratios $r:=f/g$ and $\rho:=f'/g'$, where $f(x):=\P(Z\ge x-\la/x)$ and $g(x):=e^\la\P(Z\ge x)$ for $x>0$. Then
$$\rho(x)=e^{-\frac{\la^2}{2 x^2}}(1+\tfrac\la{x^2})\quad\text{and}\quad
\rho'(x)=
\big(\la^2-(2-\la)x^2\big)\,\la\, x^{-5}\,e^{-\frac{\la^2}{2 x^2}},
   $$
so that $\rho\u\d$ on $(0,\infty)$. By part (ii) of  Proposition~\ref{prop:special-ud}, one now has  $r\d$  or $\u\d$ on $(0,\infty)$, and hence $r\d$  or $\u\d$ on $(z_V,\infty)$. 
Besides, $r(z_V)=1.020\ldots>1$ and, by l'Hospital's rule for limits, $r(\infty-)=\rho(\infty-)=1$. Therefore, $r>1$ on $(z_V,\infty)$. 
Now Lemma~\ref{lem:less} follows.
\end{proof}

\begin{proof}[Proof of Proposition~\ref{prop:W tilde}]
Note that $\la>1$ and
$\P(|Z|\ge x-\la/x)=1$ if $0<x\le\sqrt\la$. 
Therefore, $\widetilde W(x)=1$ if $0<x\le1$ and $\widetilde W(x)=\frac1{x^2}$ if $1\le x\le\sqrt\la$. 

To compare $\frac1{x^2}$ and $\P(|Z|\ge x-\la/x)$ with each other for $x\ge\sqrt\la$,
consider the ratios $r:=f/g$ and $\rho:=f'/g'$, where $f(x):=\P(|Z|\ge x-\la/x)$ and $g(x):=\frac1{x^2}$. 
Let $\vpi$ denote the density function of the standard normal distribution. Then 
$$
\rho'(x)=\Big(\tfrac{\la^3}{x^6}+\tfrac{(\la+1)
   \la}{x^4}+\tfrac{3-\la}{x^2}-1\Big)\, x^4\,
   \vpi(x-\la/x)
   $$
for $x>\sqrt\la$, so that $\rho\u\d$ on $(\sqrt\la,\infty)$. By Proposition~\ref{prop:special-ud}, this implies that $r\d$  or $\u\d$ on $(\sqrt\la,\infty)$. 
Besides, $r(\sqrt\la)=\la>1$ and $r(\infty-)=0$. 
Now Proposition~\ref{prop:W tilde} follows.
\end{proof}

\begin{remark}
It is seen from the proof of Lemma~\ref{lem:less} and expressions \eqref{eq:V} and \eqref{eq:W} for $V$ and $W$ that 
$$W(x)/V(x)\to1\quad\text{as}\quad x\to\infty.$$
\end{remark}


\end{document}